\documentclass[12pt]{amsart}
\usepackage{amscd,amssymb}
\usepackage[graph,frame,poly,arc]{xy}  
\usepackage[plainpages,backref,urlcolor=blue]{hyperref}

\usepackage{tikz}
\usetikzlibrary{calc}

\theoremstyle{plain}
\newtheorem{thm}[subsection]{Theorem}
\newtheorem{lem}[subsection]{Lemma}
\newtheorem{prop}[subsection]{Proposition}

\theoremstyle{definition}
\newtheorem{rk}[subsection]{Remark}

\numberwithin{equation}{section}
\newcommand{\al}{{\alpha}}

\newcommand{\PPP}{{\mathcal P}}

\newcommand{\Z}{\mathbb{Z}}

\newcommand{\C}{\mathbb{C}}
\newcommand{\PP}{\mathbb{P}}

\DeclareMathOperator{\dd}{d}

\DeclareMathOperator{\mult}{mult}

\begin{document}

\title [Free and nearly free curves from conic pencils]
{Free and nearly free curves from conic pencils}

\author[Alexandru Dimca]{Alexandru Dimca$^1$}
\address{Universit\'e C\^ ote d'Azur, CNRS, LJAD, France  }
\email{dimca@unice.fr}

\thanks{$^1$ This work has been supported by the French government, through the $\rm UCA^{\rm JEDI}$ Investments in the Future project managed by the National Research Agency (ANR) with the reference number ANR-15-IDEX-01.} 


\subjclass[2010]{Primary 14H50; Secondary  14B05, 13D02, 32S35, 32S40, 32S55}

\keywords{plane curves; conic pencil; free curve; syzygy; Alexander polynomial}

\begin{abstract} We construct some infinite series of free and nearly free curves using pencils of conics with a base locus of cardinality at most two. These curves have an interesting topology, e.g. a high degree Alexander polynomial that can be explicitly determined,  a Milnor fiber homotopy equivalent to a bouquet of circles, or an irreducible translated component in the characteristic variety of their complement. Monodromy eigenspaces in the first cohomology group of the corresponding Milnor fibers are also described in terms of explicit differential forms.

\end{abstract}
 
\maketitle


\section{Introduction} 

Let $S=\C[x,y,z]$ be the graded polynomial ring in the variables $x,y,z$ with complex coefficients and let $C:f=0$ be a reduced curve of degree $d$ in the complex projective plane $\PP^2$. The minimal degree of a Jacobian syzygy for $f$ is the integer $mdr(f)$
defined to be the smallest integer $r \geq 0$ such that there is a nontrivial relation
\begin{equation}
\label{rel_m}
\rho: af_x+bf_y+cf_z=0
\end{equation}
among the partial derivatives $f_x, f_y$ and $f_z$ of $f$ with coefficients $a,b,c$ in $S_r$, the vector space of  homogeneous polynomials of degree $r$. Such a curve $C$ is free (resp. nearly free) if  the graded $S$-module of Jacobian syzygies $AR(f)\subset S^3$ consisting of all relations of type \eqref{rel_m} is free (resp. has a very special minimal resolution), see  \cite{Dmax} for details.
The knowledge of the total Tjurina number of $C$, denoted by $\tau(C)$, which is the sum of the Tjurina numbers $\tau(C,p)$ for all the singular points $p$ of $C$, and of
the invariant $mdr(f)$, allows one to decide if the curve $C$ is free or nearly free. Indeed, the curve $C$  is free (resp. nearly free) if and only if $\tau(C)=(d-1)^2-r(d-r-1)$
(resp. $\tau(C)=(d-1)^2-r(d-r-1)-1$), where $r=mdr(f)$, see \cite{duPCTC,Dmax}.

Assume from now on that $C$ is not a  union of lines  passing through one point, which is equivalent to $mdr(f)>0$. When $C$ is a free (resp. nearly free) curve in the complex projective plane $\PP^2$, then the exponents of $C$, denoted by $d_1 \leq d_2$, satisfy  $d_1 =mdr(f)\geq 1$  and
one has
\begin{equation}
\label{sum1}
d_1+d_2=d-1,
\end{equation}
(resp. $d_1+d_2=d$).  For more on free hypersurfaces and free hyperplane arrangements see 
\cite{KS, OT, T, DHA, DStexpo}.

If the curve $C$ is reducible, one also calls it  a curve arrangement. When the curve $C$ can be written as the union of at least three members of a pencil of curves, we say that $C$ is a curve arrangement of pencil type. Such arrangements play a key role in the theory of line arrangements, see for instance \cite{FY, DHA}, and their relation to freeness was considered in \cite{DMich, JV}.

From the topological view-point, we consider the complement  $U= \PP^2 \setminus C$ and let $F:f=1$ be the corresponding Milnor fiber in $\C^3$, with the usual monodromy action $h:F \to F$. One can also consider  the characteristic polynomials of the monodromy, namely
\begin{equation} 
\label{Delta}
\Delta^j(t)=\det (t\cdot Id -h^j|H^j(F,\C)),
\end{equation} 
for $j=0,1,2$. It is clear that, when the curve $C$ is reduced,  one has $\Delta^0_C(t)=t-1$, and moreover
\begin{equation} 
\label{Euler}
\Delta^0_C(t)\Delta^1_C(t)^{-1}\Delta^2_C(t)=(t^d-1)^{\chi(U)},
\end{equation} 
where $\chi(U)$ denotes the Euler characteristic of the complement $U$, see for instance  \cite[Proposition 4.1.21]{D1}. Recall that
\begin{equation} 
\label{Euler1}
 \chi(U)=(d-1)(d-2)+1-\mu(C),
\end{equation}  
where  $\mu(C)$ is total Milnor number of $C$, which is the sum of the Milnor numbers $\mu(C,p)$ for all the singular points $p$ of $C$. 
 It follows that the polynomial $\Delta(t)=\Delta^1(t)$, also called the Alexander polynomial of $C$,  determines the remaining polynomial $\Delta^2(t)$. When $\chi(U) \leq 0$, a situation described for instance in \cite{GuPa, JoSt, WV} and occurring in Theorems \ref{thm1}, \ref{thm2} below, then $\Delta(t)$ is quite large.
Recall also the Hodge spectrum definition
\begin{equation} 
\label{sp1}
Sp^j(f)=\sum_{\al>0}n^j_{f,\al}t^{\al} 
\end{equation} 
for $j=0,1$, where 
$$n^j_{f,\al}=\dim Gr_F^pH^{2-j}(F,\C)_{\lambda}$$
with $F$ denoting  the Hodge filtration in $Gr_F$, $p= \lfloor 3-\al \rfloor$ and $\lambda=\exp(-2 \pi i \alpha)$, see \cite{BS,DS1,Sa1,DHA}. It is clear that $Sp^1(f)$ determines the Alexander polynomials $\Delta(t)$.

It is an interesting question to see how the freeness of a curve $C$ is reflected in the topological properties of $U$ and $F$, see for instance \cite{ArD}.
In this note we show that many interesting free and nearly curves can be obtained from pencils of conics. Our examples go beyond the papers \cite{DMich, JV}, where pencils are considered
mostly under the hypothesis that the base locus is smooth and no member in the pencil has non-isolated singularities. The freeness of conic-line arrangements is also discussed in \cite{STo}, from a different perspective and with a different aim. The topology of the complements of some conic arrangements is discussed in \cite{AGT, BT}.

Consider first the following conic pencil with one point base locus:
\begin{equation}
\label{conics1}  
\PPP_1: tx^2+s(xz+y^2)=0
\end{equation}

\begin{figure}[h]
\centering
\begin{tikzpicture}[scale=1.5]

\draw[style=thick,color=blue] (2,2) ellipse (3cm and 1cm);
\draw[style=thick,color=green] (1.5,2) ellipse (2.5cm and 0.7cm);
\draw[style=thick,color=red] (1,2) ellipse (2cm and 0.5cm);
\node at (-1.4,2){$b$};
\end{tikzpicture}
\caption{The pencil $\PPP_1$}
\label{fig:p1}

\end{figure}

In this pencil, there is a double line $2L$, where $L:x=0$, and all the other members are smooth conics, meeting just at the base point $b=(0:0:1)$. Using this pencil we construct the following curves:
\begin{equation}
\label{curves1}  
C_{2m}: f= x^{2m}+(xz+y^2)^m=0 \text{ and } C_{2m+1}: f= x(x^{2m}+(xz+y^2)^m)=0.
\end{equation}
These curves have been essentially introduced by C.T.C. Wall in \cite[Chapter 7, Section 7.5, p. 179]{CTC} (where the common tangent is $x=0$ and not $y=0$ as claimed) and independently by Arkadiusz P\l oski in \cite{Po}. These authors showed that these curves have a maximal possible Milnor number at $b$, namely
\begin{equation}
\label{Milnor1}  
\mu(C_{d},b) = (d-1)^2- \lfloor \frac{d}{2} \rfloor 
\end{equation}
in the class of all plane curves $C$ of degree $d$, with $\mult_bC <d$. Then Jaesun Shin has shown that if we consider the total Milnor number $\mu(C)$ of  plane curves $C$ of degree $d$, we get the same result, see \cite{Shin}.
Since very singular plane curves tend to be free, the first claim of 
our first main result  is not surprising. The other properties of these curves listed below are quite unusual in our opinion.

\begin{thm}
\label{thm1} 
Consider the curves $C_d$ defined in \eqref{curves1}, for $d\geq 3$. Then the following holds.
\begin{enumerate}

\item  The curves $C_d$ are free with exponents $d_1=1$ and $d_2=d-2$. In particular, the global Tjurina number $\tau(C_d)=(d-1)^2-(d-2)$ is maximal in the class of all plane curves $C$ of degree $d$, with $\mult_bC <d$.

\item  The complement $U$ satisfies $b_2(U)=0$. Moreover, $U$ is  homotopy equivalent to a bouquet of circles $\vee S^1$ if and only if $d$ is odd. In addition, the Euler characteristic $\chi(U)$ is given by 
$$\chi(U)=2-d+\lfloor \frac{d}{2} \rfloor .$$

\item When $d$ is odd, then the Milnor fiber $F$ is  homotopy equivalent to a bouquet of circles $\vee S^1$, and hence the corresponding Alexander polynomial $\Delta(t)$ of $C_d$ is given by
$$\Delta(t)=(t-1)(t^d-1)^{-\chi(U)}.$$

\item When $d=2m$ is even, the Alexander polynomial $\Delta(t)$ of $C_d$ is determined by the  Hodge spectrum
$$Sp^1(f)=\sum_{j=3,d-1}\lfloor \frac{j-1}{2} \rfloor (t^{1+\frac{j}{d}}+ t^{3-\frac{j}{d}})+(m-1)t^2  .$$

\end{enumerate}

\end{thm}

\begin{rk}
\label{rk1} (i) The   Hodge spectrum $Sp^1(f)$ in the case $d$ odd follows easily from Proposition \ref{propSp} and Lemma \ref{lem1} below.

\noindent (ii) The characteristic polynomial $\Delta^2(t)$ of $C_d$ is non-trivial when $d=2m$ is even. For instance, it follows easily from Theorem \ref{thm1} (4) and formula \eqref{Euler}, that
$\Delta^2(t)=\Phi_2(t)\Phi_6(t)$ (resp. $\Delta^2(t)=\Phi_8(t)$) for $d=6$ (resp. $d=8$).
Here $\Phi_j$ denotes the $j$-th cyclotomic polynomial.

\noindent (iii) Since $\tau(C_d) < \mu(C_d)$ for $d \geq 5$, it follows that the singularity $(C_d,b)$ is not weighted homogeneous in this range. Note also that the computation of $\tau(C_d)$ is rather difficult without using the freeness of the curve $C_d$.
\end{rk}

Consider next the following conic pencil with a two point point base locus:
\begin{equation}
\label{conics2}  
\PPP_2: txz+sy^2=0
\end{equation}

\begin{figure}[h]
\centering
\begin{tikzpicture}[scale=1]

\draw[style=thick,color=blue] (2,2) ellipse (3cm and 1cm);
\draw[style=thick,color=green] (2,2) ellipse (3cm and 1.5cm);
\draw[style=thick,color=red] (2,2) ellipse (3cm and 2cm);
\node at (-1.4,2){$b$};
\node at (5.4,2){$b'$};
\end{tikzpicture}
\caption{The pencil $\PPP_2$}
\label{fig:p2}

\end{figure}

This pencil is considered also in Shin's paper \cite{Shin} mentioned above.
In this pencil, there is a double line $2L'$, where $L':y=0$, a singular conic $Q': xz=0$,
and all the other members are smooth conics, meeting  at the base points $b=(0:0:1)$ and $b'=(1:0:0)$. Using this pencil we construct the following curves:
\begin{equation}
\label{curves2}  
C'_{2m}: f=xz[(xz)^{m-1}+y^{2m-2}] =0 \text{ and } C'_{2m+1}: f= x[(xz)^{m}+y^{2m}]=0.
\end{equation}
These curves present a dramatic change in their Alexander polynomials when we pass from an even degree to an odd one.
\begin{thm}
\label{thm2} 
Consider the curves $C'_d$ defined in \eqref{curves2}, for $d\geq 3$. Then the following holds.
\begin{enumerate}

\item  The curves $C'_d$ are free with exponents $d_1=1$ and $d_2=d-2$. In particular, the global Tjurina number $\tau(C'_d)=(d-1)^2-(d-2)$ is maximal in the class of all plane curves $C$ of degree $d$, with $\mult_bC <d$. Moreover, all the singularities of the curve $C'_d$ are weighted homogeneous and hence $\mu(C'_d)=\tau(C'_d)$.

\item  The complement $U'=\PP^2 \setminus C'_d$ satisfies $\chi(U')=0$. More precisely, one has
$$b_1(U')=\lfloor \frac{d}{2} \rfloor   \text{   and   } \  \   b_2(U')=\lfloor \frac{d}{2} \rfloor-1.$$

\item When $d=2m+1$ is odd,  one has $H^1(U',\C)=H^1(F',\C)$, where $F'$ denotes the corresponding Milnor fiber. Hence the corresponding Alexander polynomial $\Delta(t)$ of $C'_d$ is given by
$$\Delta(t)=(t-1)^{b_1(U')}.$$

\item When $d=2m$ is even, the Alexander polynomial $\Delta(t)$ of $C'_d$ is determined by the  Hodge spectrum
$$Sp^1(f)=\sum_{j=3,d-1}\lfloor \frac{j-1}{2} \rfloor (t^{1+\frac{j}{d}}+ t^{3-\frac{j}{d}})+mt^2  .$$

\end{enumerate}

\end{thm}

\begin{rk}
\label{rk1.5} As mentioned above, the curves $C_d$ realize the maximum value of the total Milnor number $\mu(C)$ in the class of curves $C$ of degree $d$, and  the curves realizing this maximum are essentially unique, as shown by Jaesun Shin in \cite{Shin} (and by Arkadiusz P\l oski in \cite{Po} for $\mu(C,b)$). If one asks the same question for the total Tjurina number $\tau(C)$, then Theorem \ref{thm1} (1) and Theorem \ref{thm2} (1) show that this unicity does no longer hold. For a discussion on maximal Tjurina numbers see also  \cite[Chapter 7, Section 7.5, pp. 178-179]{CTC}.
\end{rk}

Using the  pencil  \eqref{conics2} we construct also the following curves:
\begin{equation}
\label{curves3}  
C''_{2m}: f=(xz)^{m}+y^{2m} =0 \text{ and } C''_{2m+1}: f= y[(xz)^{m}+y^{2m}]=0.
\end{equation}
These curves are no longer free, but they are nearly free as defined for instance in \cite{Dmax}.
The next result shows that from a topological view-point, the behaviour of these two classes of curves can be very similar.
\begin{thm}
\label{thm3} 
Consider the curves $C''_d$ defined in \eqref{curves3}, for $d\geq 3$. Then the following holds.
\begin{enumerate}

\item  The curves $C''_d$ are nearly free with exponents $d_1=1$ and $d_2=d-1$. In particular, the global Tjurina number is given by $\tau(C''_d)=(d-1)^2-(d-2)-1$. Moreover, all the singularities of the curve $C''_d$ are weighted homogeneous and hence $\mu(C''_d)=\tau(C''_d)$.

\item  The complement $U''=\PP^2 \setminus C''_d$ satisfies $\chi(U'')=1$. More precisely, one has
$$b_1(U'')=  b_2(U'')=\lfloor \frac{d-1}{2} \rfloor.$$

\item The Alexander polynomial $\Delta(t)$ of $C'_d$ is determined by the  Hodge spectrum
$$Sp^1(f)=\sum_{j=3,d-1}\lfloor \frac{j-1}{2} \rfloor (t^{1+\frac{j}{d}}+ t^{3-\frac{j}{d}})+ \lfloor \frac{d-1}{2} \rfloor t^2  .$$

\end{enumerate}

\end{thm}

\begin{rk}
\label{rk2} The complements $U$, $U'$ and $U''$ come each with a surjective regular mapping
$\phi: U \to \PP^1 \setminus B$, $\phi': U' \to \PP^1 \setminus B'$
and $\phi'': U \to \PP^1 \setminus B''$ respectively. Here $B$, $B'$ and resp. $B''$ are finite sets of points in $\PP^1$ corresponding to the members of the pencil that occur in the given curve, see the next section for a precise description of $B$ for $C_{2m+1}$. Note that for the curve $C_{2m}$, the mapping $\phi$ has a multiple fiber $2L$, and for the curves $C'_d$ (resp. $C''_{2m}$) the mapping $\phi'$  (resp. $\phi''$) has a multiple fiber $2L'$. These multiple fibers create translated irreducible components in the corresponding characteristic varieties, as explained in \cite{D4,DHA}. 
Note also that the fundamental groups $\pi_1(U')$ (resp. $\pi_1(U'')$) are described in \cite{AGT} for $d=5,6$ (resp. $d=4,5$).

\end{rk}

In the final section we explain how our results on Alexander polynomials give in fact explicit de Rham cohomology classes in the Milnor monodromy eigenspaces 
$H^{1,0}(F)_{\lambda}=F^1H^1(F,\C)_{\lambda}$ of the Milnor fiber $F$, similar to the results in \cite{DStFor}.

We thank Arkadiusz P\l oski for useful discussions and a pointer to the references \cite{Shin, CTC}.

\section{Conic pencils with one point base locus} 

In this section we prove Theorem \ref{thm1} and give additional information on the free curves $C_d$. 

\medskip

{\it Proof of Theorem \ref{thm1}, claim (1).} 
When $d=2m$, then we have the following formulas
$$f_x=2mx^{2m-1}+mz(xz+y^2)^{m-1}, \ \ f_y=2my(xz+y^2)^{m-1} \text{ and } f_z=mx(xz+y^2)^{m-1}.$$
Hence $\rho_1: xf_y-2yf_z=0$ and hence $mdr(f)=1$. To show that $C_d$ is free, it is enough to show the existence of a Jacobian syzygy $\rho_2$ as in \eqref{rel_m} of degree $d-2$, which is not a multiple of the degree 1 syzygy $\rho_1$, see \cite{ST}. In order to do this, note that 
$$g= -2(xz+y^2)^{m-1}f_x+y^{2m-3}zf_y$$
is divisible by $f_z$, hence it yields the required syzygy $\rho_2$. When $d=2m+1$, the above syzygy $\rho_1$ still exists, and we follow the same idea. One has to consider the polynomial
$$g= -2mx(xz+y^2)^{m-1}f_x+y^{2m-1}f_y,$$
and note that this $g$ is again divisible by $f_z$. The claim about $\tau(C_d)$ is a consequence of 
the maximality of the total Tjurina number for free curves, see \cite{duPCTC, Dmax}.

\medskip

{\it Proof of Theorem \ref{thm1}, claim (2).}  Using the formulas \eqref{Euler1} and \eqref{Milnor1}, it follows that 
$$\chi(U)=b_0(U)-b_1(U)+b_2(U)=2-d+ \lfloor \frac{d}{2} \rfloor ,$$
as claimed. Next note that $b_0(U)=1$, while  $b_1(U)=m-1$ when $d=2m$, and $b_1(U)=m$ when $d=2m+1$. This implies $b_2(U)=0$. When $d=2m$, it follows that 
$$H_1(U, \Z)=\Z^{m-1} \oplus \Z/2\Z,$$
see for instance \cite[Proposition 4.1.3]{D1}, and hence $U$ is not a bouquet of circles. When $d=2m+1$, note that the mapping
$$\phi: U \to \PP^1 \setminus B,$$
given by $(x:y:z) \mapsto (x^2: xz+y^2)$, with $B=\{(0:1) \} \cup \{(1:-\zeta) \ : \ \zeta^m=-1\}$, is a locally trivial fibration with contractible fiber. Indeed, the fibers are smooth conics, homeomorphic to $S^2$, with the base point deleted. It follows that $U$ has the homotopy type of $\PP^1 \setminus B$, namely a bouquet of $m$ circles $S^1$.

\medskip

{\it Proof of Theorem \ref{thm1}, claim (3).} The Milnor fiber $F$ is a cyclic covering of $U$ of degree $d$. A covering of a 1-dimensional CW complex is still a 1-dimensional CW complex, hence the first claim follows. This implies that $b_2(F)=0$, and the formula for the Alexander polynomial $\Delta(t)$ follows from the formula \eqref{Euler}.

\medskip

{\it Proof of Theorem \ref{thm1}, claim (4).} First we state in down-to-earth terms some of our results in  \cite{DStproj}. For more details on this spectral sequence approach to the computation of the Milnor fiber monodromy we refer to \cite{DS1,DStFor,DStMFgen, Sa3,Sa4}.

For any reduced plane curve $C:f=0$, consider as in the Introduction the vector space $AR(f)_j$ of Jacobian syzygies of $f$ of degree $j$. We have a linear mapping $\delta_j: AR(f)_j \to S_{j-1}$ given by $(a,b,c) \mapsto a_x+b_y+c_z.$ For the following result we refer to  \cite{DStproj}, see especially Proposition 2.2 and Corollary 2.4.

\begin{prop}
\label{propSp} Let $C: f=0$ be a degree $d$ reduced plane curve. With the above notation, let $n_j= \dim \ker \delta_j$. Then the  Hodge spectrum $Sp^1(f)$ is given by the formula
$$Sp^1(f)=\sum_{j=3,d-1} n_{j-2}(t^{1+\frac{j}{d}}+ t^{3-\frac{j}{d}}) +b_1(U)t^2 .$$

\end{prop}

\proof In the notation from  \cite{DStproj} and \eqref{sp1} above, one has $\al=1+j/d$,
$n^1_{f,\al}=\dim E_{2}^{1,0}(f)_j$ for $1 \leq j \leq d$ and $n^1_{f,\al}=\dim E_{2}^{1,0}(f)_{d-(j-d)}$ for 
$d+1 \leq j \leq 2d-1$. Indeed, for the other values of $\al$, one clearly has $n^1_{f,\al}=0.$
Now it follows from the definition of $E_{2}^{1,0}(f)_k$ for $k \in [1,d]$, that
$$\dim E_{2}^{1,0}(f)_j=n_{j-2}$$
for $1 \leq j \leq d$, and 
$$\dim E_{2}^{1,0}(f)_{2d-j}=n_{2d-j-2}$$
for $d+1 \leq j \leq 2d-3$. Set $j'=2d-j$ and note that 
$$3-\frac{j}{d}=1+\frac{j'}{d}.$$
This shows that, for $1 \leq j < d$, the coefficient of $t^{3-\frac{j}{d}}$ has to be 
$n_{2d-j'-2}=n_{j-2}.$
\endproof
Now we come back to our curves $C_d$ and determine the sequence $n_j$. To start with, we have $AR(f)_0=0$, $n_0=0$ and $AR(f)_1$ is 1-dimensional, spanned by
$$\rho_1=(0,x,-2y).$$
Since $\delta_1(\rho_1)=0$, it follows that $n_1=1$. For $j$ satisfying $1<j<d-2$, the elements of $AR(f)_j$ are of the form 
$$\rho_h=h \rho_1=(0,xh,-2yh),$$
where $h \in S_{j-1}$. Note that $\delta_j(\rho_h)=0$ if and only if $xh_y-2yh_z=0$.
The following result is an easy exercise for the reader.

\begin{lem}
\label{lem1} 
Let $h \in S$ be a homogeneous polynomial of degree $e=j-1$ such that $xh_y-2yh_z=0$. If 
$e=2e_1$ is even, then $h=h_1(u,v)$, where $h_1 \in S_{e_1}$, $u=x^2$ and $v=xz+y^2$.
If $e=2e_1+1$ is odd, then $h=xh_1(u,v)$, where $h_1(u,v)$ is as above.

\end{lem}
Let now $j=2j_1$ be even. Then $e=2j_1-1=2(j_1-1)+1$ is odd, and it follows from Lemma \ref{lem1} that $n_j=e_1+1=j_1$. When $j=2j_1+1$ is odd, then $e=2j_1$ is even, and it follows from Lemma \ref{lem1} that $n_j=j_1+1$. It follows that
$$n_j=\lfloor \frac{j+1}{2} \rfloor,$$
which completes the proof of the last claim in Theorem \ref{thm1}.

\section{ Conic pencils with two point base locus }  

In this section we prove first Theorem \ref{thm2}.
\medskip

{\it Proof of Theorem \ref{thm2}, claim (1).} 
By a direct computation, we find the following degree one syzygies
$$\rho_1=(x,0,-z),$$
for $d$ even, and
$$\rho_1=(2mx,-y,-2(m+1)z)$$
for $d=2m+1$ odd. Then, in both cases we have $f_y=xg$ and $f_z=xh$ for some polynomials $g$ and $h$, which give rise to the degree $(d-2)$ syzygy
$$\rho_2=(0,h,-g).$$
The singularities of the curve $C'_d$ are located at $(1:0:0)$ and $(0:0:1)$ for $d$ odd, and the same plus an extra node at $(0:1:0)$ when $d$ is even. A simple computation shows that all these singularities are weighted homogeneous.

\medskip

{\it Proof of Theorem \ref{thm2}, claim (2).} This is obvious, using the formula for $\mu(C'_d)$ and \eqref{Euler1}.

\medskip

{\it Proof of Theorem \ref{thm2}, claim (3).} To prove this claim we use Proposition \ref{propSp} and show that $n_j=0$ for all $j$'s with $1\leq j \leq d-3$. As in the proof of Theorem \ref{thm1} (4), 
we have $AR(f)_0=0$, $n_0=0$ and $AR(f)_1$ is 1-dimensional, spanned by
$$\rho_1=(2mx,-y,-2(m+1)z).$$
Since $\delta_1(\rho_1)=-3$, it follows that $n_1=0$. For $j$ satisfying $1<j<d-2$, the elements of $AR(f)_j$ are of the form 
$$\rho_h=h \rho_1=(2mxh,-yh,-2(m+1)zh),$$
where $h \in S_{j-1}$. Note that $\delta_j(\rho_h)=0$ if and only if 
$$2mxh_x-yh_y-2(m+1)zh_z-3h=0.$$
The following result completes the proof of claim (3).

\begin{lem}
\label{lem2} 
Let $h \in S$ be a homogeneous polynomial  of degree $<2m-2$ such that 
$$2mxh_x-yh_y-2(m+1)zh_z-3h=0.$$
Then $h=0$.
\end{lem}

\proof Assume that a monomial $x^ay^bz^c$ enters into the polynomial $h$ with non-zero coefficient.
Then $a+b+c <2m-2$ and $2ma-b-2(m+1)c=3$. The last relation implies that $a>c$, say $a=c+e$ where $e\geq 1$. Then $b+2c+e <2m-2$ implies that
$$3=2ma-b-2(m+1)c>2me-(2m-2-e)=2(e-1)m+2+e \geq 3,$$
a contradiction. This shows that $h=0$.

\endproof

\medskip

{\it Proof of Theorem \ref{thm2}, claim (4).} As above, we have $AR(f)_0=0$, $n_0=0$ and $AR(f)_1$ is 1-dimensional, spanned now by
$$\rho_1=(x,0,-z).$$
Since $\delta_1(\rho_1)=0$, it follows that $n_1=1$. For $j$ satisfying $1<j<d-2$, the elements of $AR(f)_j$ are of the form 
$$\rho_h=h \rho_1=(xh,0,-zh),$$
where $h \in S_{j-1}$. Note that $\delta_j(\rho_h)=0$ if and only if $xh_x-zh_z=0$.
The following result is an easy exercise for the reader.

\begin{lem}
\label{lem3} 
Let $h \in S$ be a homogeneous polynomial of degree $e=j-1$ such that $xh_x-zh_z=0$. If 
$e=2e_1$ is even, then $h=h_1(u,v)$, where $h_1 \in S_{e_1}$, $u=xz$ and $v=y^2$.
If $e=2e_1+1$ is odd, then $h=yh_1(u,v)$, where $h_1(u,v)$ is as above.

\end{lem}
It follows as above that
$$n_j=\lfloor \frac{j+1}{2} \rfloor,$$
which completes the proof of the last claim in Theorem \ref{thm2}.

\medskip

Next we prove Theorem \ref{thm3}.

\medskip

{\it Proof of Theorem \ref{thm3}, claim (1).} 
By a direct computation, we find the  degree one syzygy
$$\rho_1=(x,0,-z).$$
Consider then the degree $d-1$ syzygies given by the Koszul relations
$$\rho_2=(f_y,-f_x,0)   \text{ and }  \rho_3=(0,f_z,-f_y) .$$
It follows from \cite[Theorem 4.1]{Dmax} that $C''_d$ is a nearly free curve with exponents $d_1=1$, $d_2=d-1$ and that
$$\tau(C''_d)=(d-1)^2-(d-2)-1=d^2-3d+2.$$
The singularities of the curve $C'_d$ are located at $(1:0:0)$ and $(0:0:1)$. A simple computation shows that all these singularities are weighted homogeneous.

\medskip

{\it Proof of Theorem \ref{thm3}, claim (2).} Obvious.

\medskip

{\it Proof of Theorem \ref{thm3}, claim (3).} The same as the proof of Theorem \ref{thm2}, claim (4).

\section{De Rham cohomology of Milnor fibers } 

In this section we give explicit bases for the eigenspaces $H^{1,0}(F)_{\lambda}=F^1H^1(F,\C)_{\lambda}$, where $F$ is one of the Milnor fibers discussed above. First we fix some notation. For a syzygy $\rho=(a,b,c) \in AR(f)_j$,  we consider the differential 2-form on $\C^3$ given by
$$\omega(\rho)=a\dd y \wedge \dd z-b\dd x \wedge \dd z+c\dd x \wedge \dd y\in \Omega^2,$$
where $\Omega^j$ denotes the space of $j$-differential forms on $\C^3$ with polynomial coefficients.
Then we have $\deg(\omega(\rho))=j+2$,  $\dd f \wedge \omega(\rho)=0$ and $\delta_j(\rho)=0$ if and only if $\dd \omega(\rho)=0$. Let $\Delta: \Omega^2 \to \Omega^1$ denote the contraction with the Euler vector field. Let $F:f=1$ be the Milnor fiber of $f$
and denote by $\iota:F \to \C^3$ the inclusion. Recall also that $h:F \to F$ is given by 
$$h(x,y,z)= \exp(2\pi i/d)\cdot (x,y,z)$$
and the monodromy action on $H^j(F,\C)$ is induced by $(h^{-1})^*$.

Consider first the curves $C_d:f=0$ from Theorem \ref{thm1}. Then take $\rho_1=(0,-x,2y)$ and set
$$\omega_1=\omega(\rho_1)=x\dd x \wedge \dd z+2y\dd x \wedge \dd y\in \Omega^2$$

If 
$e=2e_1$ is even, then let $E_e$ be the vector space of all polynomials $h=h_1(u,v)$, where $h_1 \in S_{e_1}$, $u=x^2$ and $v=xz+y^2$.
If $e=2e_1+1$ is odd, then let $E_e$ be the vector space of all polynomials  $h=xh_1(u,v)$, where $h_1(u,v)$ is as above, exactly as in Lemma \ref{lem1}.

\begin{prop}
\label{propDR1} For the Milnor fiber $F:f=1$ of the curve $C_d:f=0$, and for $k=3,...,d-1$,  the eigenspace $H^{1,0}(F)_{\lambda}$ where $\lambda=\exp (-2\pi ik/d)$ is given by the cohomology classes of  the 1-forms $\al=\iota^*(\Delta(h\omega_1))$ for $h \in E_{k-2}$.

\end{prop}

\proof In the notation from the proof of Proposition \ref{propSp}, the space
$$E_{2}^{1,0}(f)_k=E_{\infty}^{1,0}(f)_k= \{\omega \in \Omega^2 \  : \    \deg (\omega )= k, \ \  \dd f \wedge \omega= \dd \omega =0 \}$$
 is identified to $F^1H^1(F,\C)_{\lambda},$
via the map $\omega \mapsto \iota^*(\Delta (\omega))$,
see for details \cite[Remark 2.8 (i)]{DStFor} and \cite[Corollary 2.4]{DStproj}. The same proof applies for the next two similar results.
\endproof

Consider now the curves $C'_d:f=0$ from Theorem \ref{thm2} with $d$ even. Then take 
$\rho'_1=(x,0,-z)$ and set
$$\omega'_1=\omega(\rho'_1)=x\dd y \wedge \dd z-z\dd x \wedge \dd y\in \Omega^2.$$
If 
$e=2e_1$ is even, then let $E'_e$ be the vector space of all polynomials $h=h_1(u,v)$, where $h_1 \in S_{e_1}$, $u=xz$ and $v=y^2$.
If $e=2e_1+1$ is odd, then let $E'_e$ be the vector space of all polynomials  $h=yh_1(u,v)$, where $h_1(u,v)$ is as above, exactly as in Lemma \ref{lem3}.

\begin{prop}
\label{propDR2} For the Milnor fiber $F:f=1$ of the curve $C'_d:f=0$ with $d$ even, and for $k=3,...,d-1$,  the eigenspace $H^{1,0}(F)_{\lambda}$ where $\lambda=\exp (-2\pi ik/d)$ is given by the cohomology classes of  the 1-forms $\al'=\iota^*(\Delta(h\omega'_1))$ for $h \in E'_{k-2}$.

\end{prop}
Finally, by the same type of consideration, we get the following result for the curves in Theorem 
\ref{thm3}. Note that the same 1-form $\omega_1'$ and the same vector spaces $E'_k$ as above are used in this case, since the first syzygy is the same in the two situations.
\begin{prop}
\label{propDR3} For the Milnor fiber $F:f=1$ of the curve $C''_d:f=0$, and for $k=3,...,d-1$,  the eigenspace $H^{1,0}(F)_{\lambda}$ where $\lambda=\exp (-2\pi ik/d)$ is given by the cohomology classes of  the 1-forms $\al''=\iota^*(\Delta(h\omega'_1))$ for $h \in E'_{k-2}$.

\end{prop}

\end{document}